\documentclass[a4paper,12p]{article}
\usepackage{hyperref}
\title{Construction of Recurrent Fractal Interpolation Surfaces \\with Function Scaling Factors and Estimation of \\Box-counting Dimension on Rectangular Grids
}

\author{Chol-Hui Yun ,  Hui-Chol Choi ,  Hyong-Chol O\\\\
                   {\textit{Faculty of Mathematics}}, 
                        {\textit{\textbf{Kim Il Sung} University}},
                          \textit{Pyongyang, D. P. R. Korea}}

\date{Submit on 9 July, 2013}  
\usepackage{latexsym}   
\usepackage[pdftex]{graphicx}
\usepackage{amsmath}
    
\begin{document}
\maketitle      
\begin{abstract}   
We consider a construction of recurrent fractal interpolation surfaces with function vertical scaling factors and estimation of their box-counting dimension. A recurrent fractal interpolation surface (RFIS) is an attractor of a recurrent iterated function system (RIFS) which is a graph of bivariate interpolation function. For any given data set on rectangular grids, we construct general recurrent iterated function systems with function vertical scaling factors and prove the existence of bivariate functions whose graph are attractors of the above constructed RIFSs. Finally, we estimate lower and upper bounds for the box-counting dimension of the constructed RFISs. 
\\ \\ \indent
{\small\textbf{\textit{Keywords}}: Recurrent Iterated Functions System(RIFS), Fractal surface, Fractal Interpolation function(FIF), Box-counting dimension} \\ \indent
\textbf{2010 \textit{Mathematical Subject Classification}}: 37C45; 28A80;  41A05 
\end{abstract}
\section{Introduction}
\textit{Fractal interpolation surfaces} (FISs), fractal sets which are graphs of bivariate fractal interpolation functions, are being widely used in approximation theory, computer graphics, image compression, metallurgy, physics, geography, geology and so on. (See \cite{ban1, BH, fal, jac, man, YMB, zha}.) 

Barnsley \cite{ban2} defined a \textit{fractal function} as a function whose graph is an attractor of an iterated function system (IFS) and such fractal functions and fractal interpolation are widely studied to construct fractal curves and surfaces in many papers (\cite{ban1, BD1, BDD, dal, fal, GH, mal, man, mas, MY, WX, XFC, zha}). 

One {\it classical method} of construction of fractal surfaces is to use fractal curves \cite{fal, man, WX, XFC, YCO}. This approch is usefull in the case when data sets ({\it measurement data}) are given only on the {\it boundaries of the domain}. The history and recent developments of this approach are described in \cite{YCO}. 

Another ({\it direct}) method of construction of fractal surfaces is to use bivariate fractal interpolation functions. This metod is useful in the case when data sets are given on the {\it mesh-points of the whole domain}. Massopust \cite{mas} provided a construction of self-affine fractal interpolation surfaces with data set on triangular domain, where the interpolation points (data points) on the boundary are assumed coplanar. This results was generalized to allow more general boundary data and domains in \cite{GH}. Many authors studied construction methods of FISs with the data set given on rectangular grids (\cite{BD1, BDD, dal, mal, MY}). In \cite{dal} for such data sets on a rectangular grid that interpolation points on the boundary are collinear, a construction of FISs which are attractors of IFSs was provided. This was generalized in Malysz \cite{mal}, where IFS was constructed by using  constant vertical scaling factors, linear contraction transformations of domain and quadratic polynomials. In \cite{MY} the authors allowed an arbitrary data set and constructed IFS using \textit{function scaling factors} and Lipschitz transformations of domain, and estimated lower and upper bounds for the box-counting dimensions of the constructed surfaces. This estimation for box-counting dimension is improved in \cite{yun}(See \cite{FFY}, too.)

A \textit{recurrent iterated functions system} (RIFS) defined in \cite {BEH} is generalization of IFS and in \cite{BDD} they suggested a construction of \textit{recurrent fractal interpolation surfaces} (RFISs) using RIFS. This is a more flexible method of constructing fractal surfaces than using IFS and applied to the image compression (\cite{YMB}). Bouboulis and Dalla \cite{BD1} provide a general construction of \textit{recurrent fractal interpolation functions}(RFIFs) on $\textbf{R}^N$ by RIFS. In \cite{BD1, BDD} they used domain contraction transformations and {\it constant scaling} factors.

This paper is a continuation and extension of \cite{MY} and \cite{yun} where they used IFS and function vertical scaling factors. We present a flexible construction of RFISs by RIFSs with {\it function vertical scaling} factors and estimate lower and upper bounds for the box-counting dimensions of the constructed surfaces.

The remainder of the article is organized as follows: The section 2 describes construction of recurrent fractal interpolation surfaces on the rectangular grids and gives an example. In the section 3 we estimate upper and lower bounds of the box-counting dimension of RFISs constructed in section 2.   We refer to \cite{BEH, bou, YCO} for necessary preliminaries on RIFS. 
  
\section{Construction of Recurrent Fractal Interpolation Surfaces}
In this section, we construct recurrent fractal interpolation surfaces with a data set on rectangular grids. 
Let a data set on the rectangular grid be given by 
$$P=\{(x_{i},y_{j},z_{ij})\in \mathbf{R}^3;~i=0,1,\ldots,n;~j=0,1,\ldots,m \}, $$ 
$$(x_{0}<x_{1}<\ldots<x_{n};~y_{0}<y_{1}<\ldots<y_{m}). $$
Let denote $N=n\cdot m;~N_{nm}=\{1,\ldots,n\}\times \{1,\ldots,m\}$ and let
$$I_{i}=[x_{i-1},x_{i}],~J_{j}=[y_{j-1},y_{j}],~E=[x_0,x_n]\times [y_0,y_m],~E_{ij}=I_{i}\times J_{j},~(i,j)\in N_{nm}$$  
and $E_{ij}$ are called \textit{regions}.  $\mathbf{N}$ is the set of all positive integers.

Let  $l\geq 2~(l\in\mathbf{N})$. In the rectangle $E$ we choose rectangles $\tilde{E}_k(k=1,\ldots,l)$ which consist of some regions, and call  $\tilde{E}_k$ \textit{domains}. Then $\tilde{E}_k=\tilde{I}_{k}\times \tilde{J}_{k},~k=1,\ldots, l$, where $\tilde{I}_{k}$ and $\tilde{J}_{k}$ are closed intervals on the $x$ and $y$ axes, respectively. The end points of $\tilde{I}_{k}~(k\in \{1,\ldots,l\})$ coincide with some end points of intervals $I_{i}~(i=1,\ldots,n)$. If the indices of the start point and the end point of $\tilde{I}_{k}$ are respectively denoted by $s_x(k),e_x(k)$, then two mappings $s_x:\{1,\ldots,l\}\rightarrow \{1,\ldots,n\}, e_x:\{1,\ldots,l\}\rightarrow \{1,\ldots,n\}$ are well defined. For $\tilde{J}_{k}$, the two mappings $s_y:\{1,\ldots,l\}\rightarrow \{1,\ldots,m\}, e_y:\{1,\ldots,l\}\rightarrow \{1,\ldots,m\}$ are defined similarly and $\tilde{I}_{k}, \tilde{J}_{k}$ are respectively represented by 
$$\tilde{I}_{k}=[x_{s_x(k)},~x_{e_x(k)}], ~~\tilde{I}_{k}=[y_{s_y(k)},~y_{e_y(k)}].$$ 

Assume that $e_x(k)-s_x(k)\geq 2,~e_y(k)-s_y(k)\geq 2,~k=1,\cdots,l$. This means that the intervals $\tilde{I}_{k}, \tilde{J}_{k}$ include at least 2 small intervals $I_i,~J_j$ and the domain $\tilde{E}_{k}$ includes $[e_x(k)-s_x(k)]\cdot [e_y(k)-s_y(k)]$ regions.

To each region $E_{ij}$, we relate a domain $\tilde{E}_k$. This correspondence is represented by a map $\gamma :N_{nm}\rightarrow \{1,\cdots,l\}$. Throughout this paper we fix a map $\gamma$ and denote $k=\gamma(i,j)$. 

Let $L_{x,ij}:[x_{s_x(k)},x_{e_x(k)}]\rightarrow [x_{i-1}, x_i]$ and $L_{y,ij}:[y_{s_y(k)},y_{e_y(k)}]\rightarrow [y_{j-1},x_j]$, $(i,j)\in N_{nm}$ be contraction homeomorphisms. These mappings respectively map end points of $\tilde{I}_{k}, \tilde{J}_{k}$ into end points of the intervals $I_i, J_j$, that is
$$L_{x,ij}(\{x_{s_x(k)},x_{e_x(k)}\})=\{x_{i-1},x_i\},~~L_{y,ij}(\{y_{s_y(k)},y_{e_y(k)}\})=\{y_{j-1},y_j\}.$$
These mappings can be easily constructed as \cite{BDD}.

Let $s_{ij}:E_{ij}\rightarrow \mathbf {R}$ be contraction mappings on regions such that $0<|s_{ij}(x,y)|\leq s<1$, which are called \textit{vertical scaling factors}. Let $Q_{ij}:\tilde E_k\rightarrow \mathbf{R}$ be Lipschitz mappings. We define mappings $L_{ij}:\tilde{E}_k\rightarrow E_{ij}$ and $F_{ij}:\tilde{E}_k\times \mathbf{R}\rightarrow \mathbf{R}$ as folllows, respectively:
\begin{equation} \label{eq:1}   
L_{ij}(x,y)=(L_{x,ij}(x),L_{y,ij}(y)),~~~~F_{ij}(x,y,z))=s_{ij}(L_{ij}(x,y))z+Q_{ij}(x,y).
\end{equation}
\quad Now we define transformations $W_{ij}:\tilde{E}_k\times \mathbf{R}\rightarrow \mathrm E_{ij}\times \mathbf{R}(i=1,\ldots,n, j=1,\ldots,m)$ by
\begin{equation} \label{eq:2}   
W_{ij}(x,y,z)=(L_{ij}(x,y),F_{ij}(x,y,z)).
\end{equation}
Here $s_{ij}$ are taken as free unknown functions. 

For construction of RFIS, the {\it following condition} for $W_{ij}$ is {\it important}: there exists at least one continuous function $g:E\rightarrow \mathbf{R}$ interpolating the given data set $P$ such that
\begin{eqnarray}
 &&F_{ij}(x_\alpha,~y,~g(x_\alpha,~y))=g(L_{ij}(x_\alpha,~y)),~\alpha\in\{s_x(k),~e_x(k)\},~y\in [y_{s_y(k)},~y_{e_y(k)}],  \label{eq:3} \\                    
&&F_{ij}(x,~y_\beta,~g(x,~y_\beta))=g(L_{ij}(x,~y_\beta)),~x\in [x_{s_x(k)},~x_{e_x(k)}],~\beta\in\{s_y(k),~e_y(k)\}.\label{eq:4} 
\end{eqnarray}
The transformations $W_{ij}$ satisfying \eqref{eq:3} and \eqref{eq:4} can be constructed as follows, for example. \\

{\bf Example 1}\label{ex:1}. Select one Lipschitz continuous function $g_0$ interpolating the data set $P$ and let $Q_{ij}(x,~y)=g_0(L_{ij}(x,~y))-s_{ij}(L_{ij}(x,~y))\cdot g_0(x,~y),~(x,~y)\in\tilde{E}_k$, then 
\begin{eqnarray}
&&F_{ij}(x,~y,~z)=s_{ij}(L_{ij}(x,~y))\cdot (z-g_0(x,~y))+g_0(L_{ij}(x,~y)),~(x,~y)\in\tilde{E}_k.    \label{eq:5}               
\end{eqnarray}
Then the $W_{ij}$ given by \eqref{eq:2} with these $F_{ij}$ are just the needed transformations. That is why $g_0$ and the function $h$ coincided with $g_0$ in $\partial E_{ij}$ satisfy the conditions \eqref{eq:3} and \eqref{eq:4}. \\\\\indent
\textbf{Remark 1}. Then the mappings $L_{ij}$ map the vortices of the domains $\tilde{E}_k$ into the vortices of the regions $E_{ij}$. That is, for $~\alpha \in \{s_x(k),~e_x(k)\},~\beta \in \{s_y(k),~e_y(k)\},$ we have  
$$L_{ij}(x_\alpha,~y_\beta)=(x_a,~y_b)~(\mathrm{where}~a\in \{i-1,i\},~b\in \{j-1,j\}).$$
And from the conditions \eqref{eq:3} and \eqref{eq:4}, the transformations $W_{ij}$ map the data points (in $P$) given on the vortices of the domains into the data points given on the vortices of the regions, that is, for $~\alpha \in \{s_x(k),~e_x(k)\},~\beta \in \{s_y(k),~e_y(k)\},$ we have  
$$F_{ij}(x_\alpha, y_\beta, z_{\alpha\beta})=z_{ab}~(\mathrm{where}~ L_{ij}(x_\alpha,~y_\beta)=(x_a,~y_b),~a\in \{i-1,i\},~b\in \{j-1,j\}).$$
We sometimes denote $L_{ij}, W_{ij}$ by $L_{ij,k}, W_{ij,k}~(k=\gamma(i,j))$ explicitly pointing their domains. We denote Lipschitz (or contraction) constant of Lipschitz (or contraction) mapping $f$ by $L_f~(c_f)$ in what follows. \\\\\indent
We define a distance $\rho_\theta$ in $\mathbf{R}^3$ for $\theta>0$ by
$$\rho_\theta((x,y,z),(x^\prime ,y^\prime, z^\prime))=|x-x^\prime|+|y-y^\prime|+\theta|z-z^\prime|, (x,y,z),(x^\prime ,y^\prime, z^\prime)\in \mathbf{R}^3$$
as \cite{dal}. Let
$$\bar{c}_L=\max\{c_{L_{ij}}|~(i,~j)\in N_{nm}\},~\bar{L}_Q=\max \{L_{Q_{ij,k}}|~(i,~j)\in N_{nm}\}.$$
If $0<\theta<(1-\bar{c}_L)/\bar{L}_Q$, then the distance $\rho_\theta $ is equivalent to the Euclidean metric on $\mathbf{R}^3$ and $W_{ij},i=1,\ldots,n,j=1,\ldots,m$ are contraction mappings with respect to the distance $\rho_\theta$ (see \cite{BD1,BDD,dal}).

We define a row-stochastic matrix  $M=(p_{st})_{s,t=1}^N$ by 
\begin{equation}\label{c}
p_{st}=\left\{
\begin{array}{ll}
\frac{1}{a_s}, & E_{\tau ^{-1}(s)}\subset \tilde{E}_{\gamma(\tau^{-1}(t))}~,\nonumber\\
0 ,& E_{\tau ^{-1}(s)}\not\subset \tilde{E}_{\gamma(\tau^{-1}(t))}~.
\end{array}\right.
\end{equation}
Here the mapping $\tau :N_{nm}\rightarrow \{1,\ldots,N\}$ is the bijection defined by $\tau(i,j)=i+(j-1)n$ and for every fixed $s=1,\cdots,N$, the number $a_s$ indicates the number of elements of the set $\{t\in\{1,\cdots,N\}|E_{\tau ^{-1}(s)}\subset\tilde{E}_{\gamma(\tau^{-1}(t))}\}$. In other words, $a_s$ is the number of non zero elements in $s$-th row of the above row-stochastic matrix $M$. This means that $p_{st}$ is positive iff there exists a transformation $L_{ij}$ that maps $s$-th region into $t$-th region (\cite{BDD}).\\

Now we assume that $M$ is irreducible and define the {\it recurrent iterated function system (RIFS) corresponding to the given data set} $P$ by $\{\mathbf{R}^3;~M;~W_{ij},~i=1,\ldots,n,~j=1,\ldots,m\}$. Its connection matrix $C=(c_{st})_{N\times N}$ is given as follows. (Then $C$ is also irreducible.) 
 \begin{equation}\label{c}
c_{st}=\left\{
\begin{array}{ll}
1, & p_{ts}>0~,\nonumber\\
0, & p_{ts}=0~.
\end{array}\right.
\end{equation}
We denote the {\it attractor} of the RIFS constructed above by $\mathcal{A}$. The following theorem shows that $\mathcal{A}$ is a {\it recurrent fractal surface}.
\newtheorem{Thot}{Theorem}  
\begin{Thot}\label{theo1}
There exists an interpolation function $f:E\rightarrow\mathbf{R}$ of the data set $P$  whose graph is the attractor $\mathcal{A}$ of the RIFS constructed above.
\end{Thot}

{\bf (Proof)} Let  $C(E)$ be the following set:
$$C(E)=\{h\in C^0(E)|~h~\mathrm{interpolates ~the ~data ~set}~P~\mathrm{and ~satisfies ~(3), (4)}\}.$$
Then $C(E)$ is not empty set and complete metric space with respect to norm $||\cdot ||_\infty $. When $g\in C(E)$, if we define a function $Tg$ on $E$ by
$$(Tg)(x,y)=F_{ij}(L^{-1}_{ij}(x,y),~g(L^{-1}_{ij}(x,y))),~(x,y)\in E_{ij},~i=1,\cdots,n,~j=1,\cdots,m,$$
then $Tg\in C(E))$. In fact, for any fixed $(i,j)\in N_{nm}$ and for any $\alpha\in\{s_x(k),e_x(k)\})~(k=\gamma(i,j))$,
\begin{eqnarray*}
&&(Tg)(L_{ij}(x_\alpha,y))=F_{ij}(x_\alpha,y,g(x_\alpha,y))=g(L_{ij}(x_\alpha,y)),~y\in [y_{s_y(k)},~y_{e_y(k)}],\\
&&(Tg)(L_{ij}(x,y_\beta))=F_{ij}(x,y_\beta,g(x,y_\beta))=g(L_{ij}(x,y_\beta)),~x\in [x_{s_x(k)},~x_{e_x(k)}].
\end{eqnarray*}
This means that $Tg=g$ on the segments $\{(x_i,y);y\in [y_0,y_m]\}$,$\{(x,y_j);x\in [x_0,x_n]\}$, $i=0,1,\ldots ,n, j=0,1,\ldots,m$. Thus, we have  
$$ F_{ij}(x_\alpha,~y,~(Tg)(x_\alpha,~y))=F_{ij}(x_\alpha,~y,~g(x_\alpha,~y))=g(L_{ij}(x_\alpha,~y))=(Tg)(L_{ij}(x_\alpha,~y)). $$ 
Thus $Tg$ satisfies \eqref{eq:3}. Similarly, we can prove that $Tg$ satisfies \eqref{eq:4}. 

Therefore, the operator $T:C(E)\rightarrow C(E)$ is well defined and the operator $T$  is contractive from the assumption $|s_{ij}|\leq s<1$. Hence, the operator $T$  has a unique fixed point $f \in C(E)$ , which is presented by 
$$f(x,y)=F_{ij}(L^{-1}_{ij}(x,y),~f(L^{-1}_{ij}(x,y))), i=1,\ldots,n; j=1,\ldots,m. $$ 
Thus, for the graph  $Gr(f)$  of the function $f$  we obtain 
$$
Gr(f)=\bigcup_{s=1}^N \bigcup_{t \in \Lambda (s)} W_{\tau^{-1}(s)}(Gr(f|_{E_{\tau^{-1}(t)}})).
$$
Here $\Lambda (s)=\{t\in\{1,\cdots,N\}|~p_{ts}>0\}, s=1,\ldots,N$. This means that  $Gr(f)$ is the attractor of the RIFS constructed above. The uniqueness of the attractor implies $\mathcal{A}=Gr(f)$. (QED)\\

{\bf Remark 2}. Experiments often shows that If $0<|s_{ij}|<1$ outside a set with zero measure, the attractor of the RIFS constructed above becomes a recurrent fractal surface.\\

\begin{figure}[tbp]
   \centering
   \includegraphics[width=0.85\textwidth]{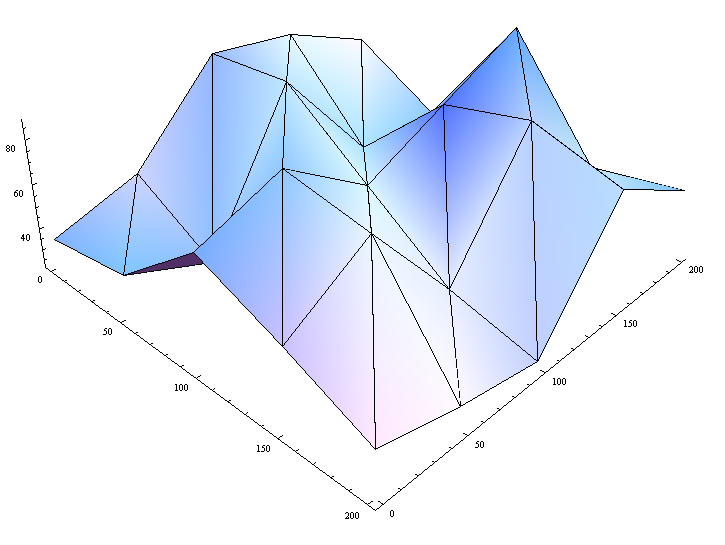}
   \caption{Data set given by table 1.}
  \label{fig1}
\end{figure}

\begin{figure}[tbp]
   \centering
   \includegraphics[width=0.6\textwidth]{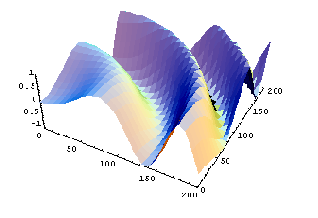}
   \caption{Vertical scaling factor $s(x, y)=\sin(10x^2+10y^2)$.}
  \label{fig2}
\end{figure}

\begin{figure}[tbp]
   \centering
   \includegraphics[width=0.85\textwidth]{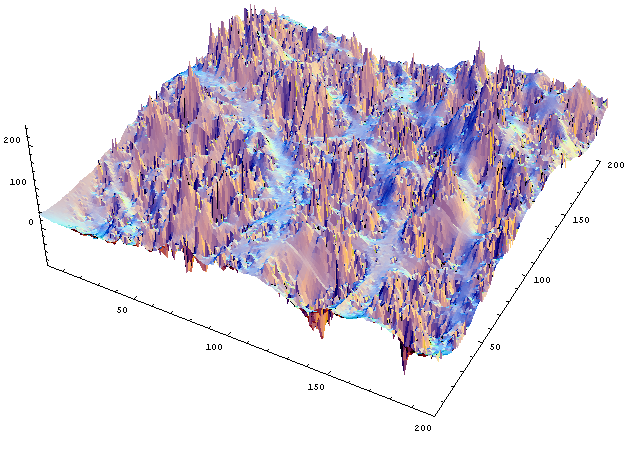}
   \caption{Recurrent Fractal Interpolation Surface constructed from the dataset in table 1.}
  \label{fig3}
\end{figure}

{\bf Example 2}. A data set is given by the following table 1, the graph of which is shown in the figure \ref{fig1}. Let $s_{ij}(x,y)=\sin(10x^2+10y^2)$ (see the figure \ref{fig2}). Let $g_0(x, y)$ be the {\it Lagrangean interpolation} function. Then the attractor of the RIFS constructed in the example 1 on page \pageref{ex:1} by \eqref{eq:2} and \eqref{eq:5} is drawn in the figure \ref{fig3}.\\
\begin{eqnarray*}
\begin{tabular}{|l|l|l|l|l|l|}
\hline
\qquad X & ~~~~0~~~~ & ~~~50~~~~ & ~~100~~~~ & ~~150~~~~ & ~~200~~~~  \\
~Y\quad &   &  &  &  &  \\
\hline
~~0~ & ~~~35~ & ~~~42~ & ~~~76~ & ~~~61~ & ~~~44~  \\
\hline
~50 & ~~~43~ & ~~~28~ & ~~~88~ & ~~~83~ & ~~~33~  \\
\hline
100 & ~~~78~ & ~~~84~ & ~~~58~ & ~~~33~ & ~~~25~  \\
\hline
150 & ~~~68~ & ~~~33~ & ~~~73~ & ~~~86~ & ~~~77~  \\
\hline
200 & ~~~47~ & ~~~29~ & ~~~88~ & ~~~43~ & ~~~54~  \\
\hline 
\end{tabular}
\end{eqnarray*}
$$\mathrm{Table 1.~Data~set~in~example~1.}$$ 


 \section{Box-counting Dimensions of Recurrent Fractal Surfaces}
In this section we estimate upper and lower bounds of the box-counting dimension of the attractor $\mathcal{A}$ of RFIS constructed in the previous section.

There exists a bi-Lipschitz homeomorphism which maps a rectangle $[0,1]\times [0,t]~(t>0)$ to any rectangle in $\mathbf{R}^2$ and box-counting dimension is invariant under bi-Lipschitz homeomorphisms. Thus we can assume that $E=[0,~1]\times [0,~m/n]$ and the end points of the regions and domains satisfy the following conditions.
\begin{eqnarray*}
&& x_{i+1}-x_i=y_{j+1}-y_j=\frac{1}{n},\quad x_{e_x(k)}-x_{s_x(k)}=y_{e_y(k)}-y_{s_y(k)}=\frac{a}{n},\\
&& i=0,1,\ldots,n-1,~j=0,1,\ldots,m-1, ~a\in\mathbf{N}, ~k=1,\ldots,l.
\end{eqnarray*}
Then there are exact $a^2$ regions in every domain. 

For $r(>0)$, we define a set $B$ of cubes as follows:
$$ \mathcal{B}_r=\left\{\left[\frac{u-1}{a^r},~\frac{u}{a^r}\right]\times \left[\frac{v-1}{a^r},~\frac{v}{a^r}\right]\times \left[b,~b+\frac{1}{a^r}\right]:~u,~v\in\mathbf{N},~b\in \mathbf{R}\right\}. $$
Let denote the smallest number of cubes in $ \mathcal{B}_r$ necessary to cover $\mathcal{A}$ by $N(\frac{1}{a^r})$  and the smallest number of $\frac{1}{a^r}$-mesh cubes that cover $\mathcal{A}$ by $N^\prime(\frac{1}{a^r})$. We can easily see that
$$ N^\prime \left(\frac{1}{a^r}\right)\leq  N\left(\frac{1}{a^r}\right)\leq 8\cdot N^\prime \left(\frac{1}{a^r}\right),$$ 
which allows us to use $N\left(\frac{1}{a^r}\right)$ to estimate the box-counting dimension of $\mathcal{A}$.

For a set $D\subset \mathbf{R}^2$, we define the {\it maximum variation} of a function $f$ on $D$ as follows:
$$R_f[D]=sup\left\{|f(x_2,~y_2)-f(x_1,~y_1)|:~(x_1,~y_1),~(x_2,~y_2)\in D\right\}.$$
\newtheorem{lem}{Lemma}
\begin{lem}\label{lem:1}                      
 Let $D$ be a rectangle in $\mathbf{R}^2$ and $W:D\times \mathbf{R}\rightarrow D\times \mathbf{R}$ the transformation of the form
\begin{equation*}
W\left( 
\begin{array}{ll}
x\\  y\\ z
\end{array}
\right)=\left(
\begin{array}{ll}
L(x,~y)\\  F(x,~y,~z)
\end{array}
\right)=
\left(
\begin{array}{ll}
L(x,~y)\\ s(L(x,~y))z+Q(x,~y)
\end{array}
\right).
\end{equation*}
Here $Q$ is Lipschitz function with the Lipschitz constants $L_Q$, $L$ is the domain contraction transformation (defined just like $L_{ij}$ in the section 2) with contraction factor $c_L$ and $s(x,~y)$ is a contraction function with $\left|s(x,~y)\right|<1$. Then for any continuous function  $f:D\rightarrow \mathbf{R}$, we have
 \begin{displaymath}
 R_{F(L^{-1},~f\circ L^{-1})}[L(D)]\leq\bar s R_f[D]+\mathrm{diam}(D)(c_s\bar f+L_Q).
 \end{displaymath}
Here $\mathrm{diam}(D)$ is a diameter of the set $D,~\bar s=\max_D\left|s(x,y)\right|, c_s$ is a contraction factor of $s(x,~y)$, $\bar f=\max_D\left|f(x,~y)\right|$.
\end{lem}

\textbf{Proof.} For $(x,~y),~(x^\prime,~y^\prime)(\in L(D))$, let denote $(\tilde x,~\tilde y)=L^{-1}(x,~y),~ (\tilde x^\prime,~\tilde y^\prime)=L^{-1}(x^\prime,~y^\prime)(\in D)$. Then we have
\begin{align*}
|F&(L^{-1}, ~f\circ L^{-1})(x,~y)-F(L^{-1},~ f\circ L^{-1})(x^\prime,~y^\prime)|= \\
&\quad=|F(L^{-1}(x,~y), ~f\circ L^{-1}(x,~y))-F(L^{-1}(x^\prime,~y^\prime),~ f\circ L^{-1}(x^\prime,~y^\prime))| \\
&\quad=|s(x,~y)f(\tilde x, ~\tilde y)+Q(\tilde x,~ \tilde y)-s(x^\prime,~y^\prime)f(\tilde x^\prime, ~\tilde y^\prime)-Q(\tilde x^\prime,~ \tilde y^\prime)| \\
&\quad=|s(x,y)f(\tilde x, \tilde y)-s(x,y)f(\tilde x^\prime, \tilde y^\prime)+s(x,y)f(\tilde x^\prime,\tilde y^\prime)-s(x^\prime,y^\prime)f(\tilde x^\prime,\tilde y^\prime)+Q(\tilde x, \tilde y)-Q(\tilde x^\prime,\tilde y^\prime)| \\
&\quad\leq\bar s R_f[D]+c_sd((x,~y),~(x^\prime,~y^\prime))\bar f+L_Q d((\tilde x,~ \tilde y),~(\tilde x^\prime, ~\tilde y^\prime)) \\
&\quad\leq\bar{s} R_f[D]+\mathrm{diam}(D)(c_s\bar f+L_Q).\qquad\qquad\qquad\qquad \mathrm{(QED)}  
\end{align*}
\quad For $N\times N$ matrix $U=(u_{ij}),~V=(v_{ij})$ we define the {\it relation} $''<~''$ by 
\begin{displaymath}
U<V\Longleftrightarrow u_{ij}<v_{ij},~i,j=1,2,\ldots ,\mathrm{N}
\end{displaymath} 
Points of a set $B$ in $\mathbf{R}^3$ are said to be $x$ (or $y$)-{\it collinear} if all the points of the set $B$ with the same $x$ (or $y$) coordinates lies on one lline. 

\begin{Thot}\label{theo2}
Let the function $f:\mathrm{E}\rightarrow\mathbf{R}$ be the interpolation function constructed in Theorem \ref{theo1}. Let $\bar{S}$ and $\underline{S}$ be $N\times N$ diagonal matrices $$
\bar{S}=\mathrm{diag}(\bar{s}_{\tau^{-1}(1)},\ldots,\bar{s}_{\tau^{-1}(N)}),~~ \underline{S}=\mathrm{diag}(\underline{s}_{\tau^{-1}(1)},\ldots, \underline{s}_{\tau^{-1}(N)}),
$$
where $\bar{s}_{\tau^{-1}(t)}=\bar s_{ij}=\max_{E_{ij}}|s_{ij}(x,~y)|,~\underline{s}_{\tau^{-1}(t)}=\underline{s}_{ij}=\min_{E_{ij}}|s_{ij}(x,~y)|,~t\in \{1,\ldots,N\}, t=\tau(i,j)$. If there exists a domain $\tilde{E}_{k_0}$ such that the interpolation points of $P\cap (\tilde{E}_{k_0}\times\mathbf{R})$ are not $x$-collinear or not $y$-collinear, then the box-counting dimension $\mathrm{dim}_B\mathcal{A}$ of the attractor $\mathcal{A}$ is estimated as follows:

1)	If $\underline \lambda>a$, then 
$$1+\log_a{\underline \lambda}\leq \mathrm{dim}_B\mathcal{A}\leq 1+\log_a{\bar{\lambda}}.$$

2)	If $\bar{\lambda}\leq a$, then 
$$\mathrm{dim}_B\mathcal{A}=2.$$
Here $\underline \lambda=\rho(\b{S}C)$ and $\bar \lambda=\rho(\bar SC)$ are spectral radii of the irreducible matrices $\underline{S}C$ and $\bar SC$, respectively.
\end{Thot}

{\bf (Proof)}. Proof of (1). We simply denote the maximum variance $R_f[\tilde{E}_{k=\gamma(i,j)}]$ by $R_{ij}$. Let denote $\frac{1}{a^r}$ by $\varepsilon_r$. Then $r\rightarrow\infty \Leftrightarrow \varepsilon _r\rightarrow 0$.

After applying once each $W_{ij}=W_{ij,~k}~(k=\gamma(i,j))$ to the interpolation points in the domain $\tilde{E}_k$, we have $(a+1)^2$ new image points of interpolation points in the region $E_{ij}$. According to the hypothesis, the interpolation points lying inside the domain $\tilde{E}_{k_0}$ are not $x$-collinear or not $y$-collinear and the $(a+1)^2$ image points in the region $E_{i_0j_0}~(k_0=\gamma(i_0,j_0))$ are not $x$-collinear or not $y$-collinear. On the other hand the connection matrix $C$ is irreducible and thus the region $E_{i_0j_0}$ is mapped into arbitrary regions $E_{ij}$ by applying the appropriately selected transformations from $\{W_{ij}:~i=1,\cdots,n;~j=1,\cdots,m\}$ several times. So in each region $E_{ij}$ there exist the $(a+1)^2$ image points of interpolation points which are not $x$-collinear or not $y$-collinear. Therefore, in each region $E_{ij}$ there are at least 3 image points of interpolation points which are not colinear and the maximum vertical distance computed only with respect to the $z$-axis from one of the 3 points to the line through other 2 points is greater than 0 (\cite{BDD}). The maximum value is called a {\it height} and denote by $H_{ij}$.

On the other hand, by Lemma \ref{lem:1}, on each region $E_{ij}$ we have
$$R_f[E_{ij}]\leq \bar{s}_{ij}R_f[\tilde{E}_{\gamma(i,j)}]+\frac{a}{n}b.$$ 
where $b=2^{\frac{1}{2}}(c_s\bar f+L_Q)$.

We define non negative vectors $\mathbf{h}_1, \mathbf{r}, \mathbf{u}_1$ and $\mathbf{i}$  as follows:
 \begin{equation*}                    
 \mathbf{h_1}=\left(
\begin{array}{ll}
H_{\tau^{-1}(1)}\\
~~~~\vdots \\
H_{\tau^{-1}(\mathrm{N})} 
\end{array} \right),~~\mathbf{r}=\left(
\begin{array}{ll}
\bar s_{\tau^{-1}(1)}R_{\tau^{-1}(1)}\\
~~~~~~\vdots \\
\bar s_{\tau^{-1}(\mathrm{N})}R_{\tau^{-1}(\mathrm{N})} 
\end{array} \right),~~\mathbf{i}=\left(
\begin{array}{ll}
1\\  \vdots \\  1 
\end{array} \right),~~\mathbf{u_1}=\mathbf{r}+\frac{a}{n}b\mathbf{i}.
\end{equation*}
Since  $\mathcal{A}$ is the graph of a continuous function defined on $E$, the smallest number of cubes in $\mathcal{B}_r$ necessary to cover $(E_{ij}\times\mathbf{R})\cap\mathcal{A}$ is greater than the smallest number of cubes in $\mathcal{B}_r$ necessary to cover vertical line with the length $H_{ij}$ and less than the smallest number of cubes in $\mathcal{B}_r$ necessary to cover the rectangular parallelepiped $E_{ij}\times[\underline{f}_{ij},~\bar{f}_{ij}]$, where
$$\underline{f}_{ij}=\min_{E_{ij}}f(x,~y),~\bar f_{ij}=\max_{E_{ij}}f(x,~y).$$
Therefore (in the bellow $[d]$ is the integer part of $d\in\mathbf{R}$),  
\begin{eqnarray*} 
&&\sum^n_{i=1}\sum^m_{j=1}[H_{ij}\varepsilon^{-1}_r]\leq N(\varepsilon_r)\leq \sum^n_{i=1}\sum^m_{j=1}\left(\left[\left(\bar s_{ij}R_{ij}+\frac{a}{n}b\right)\varepsilon^{-1}_r\right]+1\right)\left(\left[\frac{\varepsilon^{-1}_r}{n}\right]+1\right)^2, \\
&&\sum_{t=1}^N(H_{\tau^{-1}(t)}\varepsilon^{-1}_r)-N\leq N(\varepsilon_r)\leq \sum_{t=1}^N\left(\left(\bar s_{\tau^{-1}(t)}R_{\tau^{-1}(t)}+\frac{a}{n}b\right)\varepsilon^{-1}_r+1\right)\left(\left[\frac{\varepsilon^{-1}_r}{n}\right]+1\right)^2
\end{eqnarray*}
and thus if we denote $\Phi(\mathbf{a})=a_1+\dots+a_N$ for $\mathbf{a}=(a_1,\ldots,a_N)$, then we have     
$$\Phi(\mathbf{h}_1\varepsilon^{-1}_r)-N\leq N(\varepsilon_r)\leq \Phi(\mathbf{u}_1\varepsilon^{-1}_r+\mathbf{i})\left(\left[\frac{\varepsilon^{-1}_r}{n}\right]+1\right)^2,$$
where $r$ is selected so large that $\frac{1}{\mathrm{N}}>\varepsilon_r $.

After applying $W_{ij}$ twice, in each region $E_{ij}$ we have $a^2$ new small squares of side $\frac{1}{an}$, which are mapped by the transformation $W_{ij}$  from the regions $E_{i^\prime,j^\prime}$ lying inside the domain $\tilde E_{k}=\tilde E_{\gamma(i,j)}$. And since segments parallel to $z$-axis are mapped to those parallel to $z$-axis, for each region $\mathrm{E}_{ij}$  the height on these new small squares is not less than $\underline{s}_{ij}\cdot H$, where $H$ is the height on the original region $E_{i^\prime,j^\prime}$ contained in domain $\tilde{E}_k$. Therefore, the sum of maximum variances of $f$ on $a^2$ small squares of side $\frac{1}{an}$ contained in the region $\mathrm{E}_{ij}$ is not greater than $\tau(i,j)$-th coordinate of the vector 
$\mathbf{u}_2=\bar{S}C\mathbf{u}_1+\frac{a^2}{n}b\mathbf{i}$ and the sum of the heights is not less than $\tau(i,j)$-th coordinate of the vector $ \mathbf{h}_2=\underline{S}C\mathbf{h}_1$. So we have
$$
{\Phi}(\mathbf{h}_2\varepsilon^{-1}_r)-a^2N\leq N(\varepsilon_r)\leq {\Phi}(\mathbf{u}_2\varepsilon^{-1}_r+a^2\mathbf{i})\left(\left[\frac{\varepsilon^{-1}_r}{an}\right]+1\right)^2
$$
where $\frac{1}{an}>\varepsilon_r$. 

By induction we get the following conclusion: if we take $k$ such that
$$ a\varepsilon_r \geq \frac{1}{a^{k-1}n}\geq \varepsilon_r \Longleftrightarrow r-\log_an+1>k\geq r-\log_an$$
and apply $k$ times the transformations $\{W_{ij}:~i=1,\cdots,n;~j=1,\cdots,m\}$, then we get  $a^{2(k-1)}$ small squares of side $\frac{1}{a^{k-1}n}$ contained in each region $E_{ij}$ and
\begin{eqnarray} \label{eq:6}   
{\Phi}(\mathbf{h}_k\varepsilon^{-1} _r)-a^{2(k-1)}N \leq N(\varepsilon _r)\leq \mathbf{\Phi}(\mathbf{u}_k\varepsilon^{-1} _r+a^{2(k-1)}\mathbf{i})\left(\left[\frac{\varepsilon^{-1} _r}{a^{k-1}n}\right]+1\right)^2
\end{eqnarray}
where 
$$\mathbf{u}_k=\bar SC\mathbf{u}_{k-1}+\frac{a^k}{n}b\mathbf{i},~\quad \mathbf{h}_k=\underline{S}C\mathbf{h}_{k-1}.$$
Then we have
\begin{eqnarray*}		
&& \mathbf{u}_k=(\bar SC)^{k-1}\mathbf{r}+(\bar SC)^{k-1}\frac{a}{n}b\mathbf{i}+(\bar SC)^{k-2}\frac{a^2}{n}b\mathbf{i}+\ldots+(\bar SC)\frac{a^{k-1}}{n}b\mathbf{i}+\frac{a^k}{n}b\mathbf{i},\\
&& \mathbf{h}_k=(\underline{S}C)^{(k-1)}\mathbf{h}_1.	
\end{eqnarray*}
Since $\underline{S}C$, $\bar{S}C$ are non-negative irreducible matrix, from Frobenius's theorem (see \cite{bou, YCO}) there are strictly positive eigenvectors $\underline{\mathbf{e}}$, $\bar{\mathbf{e}}$ of $\underline{S}C$,  $\bar{S}C$ which correspond to eigenvalues $\underline{\lambda} =\rho(\underline{S}C)$, $\bar{\lambda}=\rho(\bar{S}C)$ of $\underline{S}C$, $\bar{S}C$ and we can choose   $\bar{\mathbf{e}}$ , $\underline{\mathbf{e}}$ so that 
 \begin{eqnarray*}
0<\mathbf{\underline{\mathbf{e}}}<\mathbf{h}_1,\quad \mathbf{r}\leq \mathbf{\bar e}, \quad b\mathbf{i}<n\mathbf{\bar e}
\end{eqnarray*}
Then by \eqref{eq:6}, we have
\begin{align}
N(\varepsilon _r)&\leq {\Phi}\left(\mathbf{u}_k\varepsilon^{-1} _r+a^{2(k-1)}\mathbf{i}\right)\left(\left[\frac{\varepsilon^{-1} _r}{a^{k-1}n}\right]+1\right)^2  \nonumber \\
&\leq {\Phi}\left(\mathbf{u}_k\varepsilon^{-1} _r+a^{2(k-1)}\mathbf{i}\right)(a+1)^2 \nonumber\\
&\leq {\Phi}\left( (\bar S\mathrm{C})^{k-1}\mathbf{r}\varepsilon^{-1} _r+(\bar S\mathrm{C})^{k-1}\frac{a}{n}b\mathbf{i}\varepsilon^{-1} _r+(\bar S\mathrm{C})^{k-2}\frac{a^2}{n}b\mathbf{i}\varepsilon^{-1} _r+\ldots \right.\nonumber\\
&\quad \left.+(\bar S\mathrm{C})\frac{a^{k-1}}{n}b\mathbf{i}\varepsilon^{-1} _r+\frac{a^k}{n}b\mathbf{i}\varepsilon^{-1} _r+a^{2(k-1)}\mathbf{i}\right)(a+1)^2  \nonumber\\
&\leq {\Phi}\left( (\bar S\mathrm{C})^{k-1}\mathbf{\bar e}\varepsilon^{-1} _r+(\bar S\mathrm{C})^{k-1}\mathbf{\bar e}a\varepsilon^{-1} _r+(\bar S\mathrm{C})^{k-2}\mathbf{\bar e}a^2\varepsilon^{-1}_r+\ldots \right.\nonumber\\
&\quad \left. +(\bar S\mathrm{C})\mathbf{\bar e}a^{k-1}\varepsilon^{-1} _r+\mathbf{\bar e}a^k\varepsilon^{-1} _r+a^{2(k-1)}\mathbf{i}\right)(a+1)^2  \nonumber\\
&=\left\{\bar{\lambda}^{k-1}{\Phi}(\mathbf{\bar e})\varepsilon^{-1}_r+\bar {\lambda}^{k-1}{\Phi}(\mathbf{\bar e})a\varepsilon^{-1}_r+\bar \lambda^{k-2}{\Phi}(\mathbf{\bar e})a^2\varepsilon^{-1}_r+\ldots \right.\nonumber\\
&\quad \left.+\bar {\lambda}{\Phi} (\mathbf{\bar e})a^{k-1}\varepsilon^{-1} _r+{\Phi} (\mathbf{\bar e})a^k\varepsilon^{-1} _r+a^{2(k-1)}{\Phi}(\mathbf{i})\right\}(a+1)^2  \nonumber\\
&\leq \left\{\bar {\lambda}^{r-\nu }\bar \mu \varepsilon^{-1} _r+\bar {\lambda }^{r-\nu}a\bar \mu \varepsilon^{-1} _r+\bar {\lambda }^{r-\nu -1}a^2\bar \mu \varepsilon^{-1}_r+\ldots \right.\nonumber\\
&\left. \quad+\bar {\lambda}a^{r-\nu} \bar \mu \varepsilon^{-1} _r+ a^{r-\nu+1}\bar \mu \varepsilon^{-1} _r+a^{2(r-\nu)}N\right\}(a+1)^2.    \label{eq:7} 
\end{align}
where  $\nu=\log_an, ~\bar \mu ={\Phi}(\mathbf{\bar e})$.

On the other hands, since $(\underline{S}C)_{ij}\leq (\bar{S}C)_{ij}$ for $(i,j)(\in N_{nm}$  , from Frobenius's theorem we have $\underline{\lambda}\leq \bar{\lambda}$. If $\underline{\lambda}>a$, then $1>\frac{a}{\underline{\lambda}}\geq \frac{a}{\bar{\lambda}}$ and thus we obtain
 \begin{eqnarray*} 
 N(\varepsilon_r)&\leq& \bar{\lambda}^{r-\nu}\bar\mu\varepsilon^{-1}_r \left(1+a+\frac{a^2}{\bar{\lambda}}+\ldots+\frac{a^{r-\nu +1}}{\bar{\lambda}^{r-\nu}}+\frac{a^{r-2\nu}N}{\bar{\lambda}^{r-\nu}\bar\mu}\right)(a+1)^2 \\
&=&\bar{\lambda}^{r}\varepsilon^{-1}_r\bar{\lambda}^{-\nu}\bar\mu \left(1+a+\frac{1-\left({a}/{\bar{\lambda}}\right)^{r-\nu +1}}{1-{a}/{\bar{\lambda}}}+\frac{a^{r-2\nu}N}{\bar {\lambda}^{r-\nu}\bar \mu }\right)(a+1)^2.
\end{eqnarray*}
 Let denote
 \begin{eqnarray*}						
\delta(r)=\bar{\lambda}^{-\nu}\bar\mu \left(1+a+\frac{1-\left({a}/{\bar{\lambda}}\right)^{r-\nu +1}}{1-{a}/{\bar{\lambda}}}+\frac{a^{r-2\nu}N}{\bar {\lambda}^{r-\nu}\bar \mu }\right)(a+1)^2
\end{eqnarray*}
Then $\delta (r)>0$ and $\frac{\log N(\varepsilon_r)}{-\log\varepsilon_r} \leq 1+\log_a\bar{\lambda}+\frac{1}{r}\log_a\delta(r)$, thus we have
\begin{eqnarray}\label{eq:8}  
{dim}_B\mathcal{A}=\lim_{\varepsilon_r\rightarrow 0}\frac{\log N(\varepsilon_r)}{-\log\varepsilon_r}\leq 1+\log_a{\bar {\lambda}}.	
\end{eqnarray}	

By \eqref{eq:6}, we have
\begin{eqnarray*}
 N(\varepsilon_r)&\geq& {\Phi}(\mathbf{h}_k\varepsilon^{-1}_r)-a^{2(k-1)}N~=~{\Phi}((\underline{S}C)^{k-1}\mathbf{h}_1\varepsilon^{-1} _r)-a^{2(k-1)}N  \\
& \geq& {\Phi}((\underline{S}C)^{k-1} \underline{\mathbf{e}}\varepsilon^{-1}_r)-a^{2(k-1)}N~=~\underline{\lambda}^{k-1}\underline{\mu}\varepsilon^{-1}_r-a^{2(k-1)}N  \\
& \geq& \underline{\lambda}^{r-\nu-1}\underline{\mu}\varepsilon^{-1}_r-a^{r-2\nu}N\varepsilon^{-1}_r  \\ 
&=& \varepsilon^{-1}_r\underline{\lambda}^{r}\underline{\lambda}^{-\nu-1} \left(\underline{\mu}-\frac{a^{r-2\nu}N}{\underline{\lambda}^{r-\nu-1}}\right), 
\end{eqnarray*}
where $\underline{\mu}={\Phi}(\mathbf{\underline{\mathbf{e}}})$. Since $\underline{\lambda}>a$, there is $r_0$ such that 
$$\eta (r):=\underline{\lambda}^{-\nu-1}\left(\underline{\mu}-\frac{a^{r-2\nu}{N}}{\underline{\lambda}^{r-\nu-1}}\right),\quad \mathrm{for~any~} r>r_0. $$ 
Therefore for $r>r_0$ we have $\frac{\log N(\varepsilon _r)}{-\log\varepsilon _r}\geq 1+\log_a\underline{\lambda}+\frac{1}{r}\log_a{\eta (r)}$ and thus
\begin{eqnarray}\label{eq:9}  
{\dim}_B\mathcal{A}=\lim_{\varepsilon _r\rightarrow 0}\frac{\log N(\varepsilon _r)}{-\log\varepsilon _r}\geq 1+\log_a{\b{$\lambda$}}.	
\end{eqnarray} 
 By \eqref{eq:8} and \eqref{eq:9}, if $\underline{\lambda}>a$, then we get
$$1+\log_a{\underline \lambda}\leq \mathrm{dim}_B\mathcal{A}\leq 1+\log_a{\bar{\lambda}}.$$
\qquad Proof of 2). If $\bar{\lambda}\leq a$ , then by \eqref{eq:7}
\begin{align*}	
N(\varepsilon _r)&\leq \{\bar{\lambda}^{r-\nu}\bar \mu \varepsilon^{-1}_r+\bar{\lambda}^{r-\nu}a\bar \mu \varepsilon^{-1}_r+\bar{\lambda}^{r-\nu-1}a^2\bar \mu \varepsilon^{-1}_r+\ldots +\bar{\lambda}a^{r-\nu}\bar \mu \varepsilon^{-1}_r+\nonumber \\
&\quad+a^{r-\nu +1}\bar \mu \varepsilon^{-1}_r+a^{2(r-\nu )}{N}\}(a+1)^2 \nonumber \\
&\leq \{a^{r-\nu}\bar \mu \varepsilon^{-1}_r+a^{r-\nu }\bar \mu \varepsilon^{-1}_r+a^{r-\nu +1}\bar \mu \varepsilon^{-1}_r+\ldots +a^{r-\nu +1}\bar \mu \varepsilon^{-1}_r+\nonumber \\
&\quad+a^{r-\nu +1}\bar \mu \varepsilon^{-1}_r+a^{r-2\nu }{N}\varepsilon^{-1}_r\}(a+1)^2 \nonumber \\
&\leq \varepsilon^{-2}_r\{a^{-\nu}\bar\mu +(r-\nu +1)a^{-\nu +1}\bar\mu +a^{-2\nu} {N}\}(a+1)^2,
\end{align*}	
Hence, we have
\begin{displaymath}
{\dim}_B\mathcal{A}=\lim_{\varepsilon _r\rightarrow 0}\frac{\log N(\varepsilon _r)}{-\log\varepsilon _r}\leq 2+\frac{1}{r}\log_a\{(a^{-\nu}\bar\mu +(r-\nu +1)a^{-\nu +1}\bar\mu +a^{-2\nu} {N}\}(a+1)^2.
\end{displaymath}
On the other hands, since $\mathcal{A}$ is the surface in $\mathbf{R}^3$, we have ${\dim}_B\mathcal{A} \geq 2$. Hence ${\dim}_B\mathcal{A}=2 $. (QED) \\

\textbf{Remark 3}. In the case where $s_{ij}(x,~y)=s_{ij}$(constant), if $\bar{\lambda}= \b{$\lambda$}>a$ , then ${\dim}_B\mathcal{A}=1+\log_a{\lambda}$. This is the estimation of Box-counting dimension of RFISs in \cite{BDD}.

\end{document}